\documentclass[11pt]{article}
\usepackage{amsfonts,amssymb}

\def\demo{\underline{Proof:}}
\def\finde{\hfill\fbox{\mbox {}} \\}

\newtheorem{df}{Definition}
\newtheorem{Th}[df]{Theorem}
\newtheorem{prop}[df]{Proposition}
\newtheorem{lema}[df]{Lemma}

\begin{document}
\author{Daniel Sadornil}

\begin{center}
{\Large{Elliptic curves with rational subgroups of order three}}
\end{center}
\begin{center} D. Sadornil\footnote{Departamento of
Matem\'aticas, Universidad de Salamanca, Plaza de la Merced 1-4,
37008 Salamanca, Spain. {\small email: {\tt sadornil@usal.es}}\\
Supported by MTM2004-00876} \\
\vspace*{1cm} Version In.1 \\

\today
\end{center}
\begin{abstract}

In this article we present a characterization of elliptic curves
defined over a finite field $\mathbb{F}_q$ which possess a rational
subgroup of order three. There are two posible cases depending on
the rationality of the points in these groups. We show that for
finite fields $\mathbb{F}_q$, $q \equiv -1 \bmod{3}$, all elliptic
curves with a point of order 3, they have another rational subgroup
whose points are not defined over the finite field. If $q \equiv 1
\bmod{3}$, this is no true; but there exits a one to one
correspondence between curves with points of order 3 and curves with
rational subgroups whose points are not rational.

\end{abstract}

\section{Introduction}

Elliptic curves over finite fields have been largely studied in the
literature (for example, see \cite{Men:ell} and the references cited
there). For an elliptic curve, we have an addition law which endows
to the rational points (those defined over the finite field)  a
group structure. This group is an abelian group of rank 1 or 2. The
type of the group is $(n_1,n_2)$, i.e. $E(\mathbb{F}_q) \cong
\mathbb{Z}/n_1\mathbb{Z} \times \mathbb{Z}/n_2\mathbb{Z}$ where $n_2
\mid n_1$ and furthermore $n_2 \mid q-1$.\\

This structure gives us some information about the torsion points of
an elliptic curve, but it is not easy to find them. This problem
arises since the computation of  the cardinal of an elliptic curve
is, in general, very hard. Nevertheless, we have an easy method to
test if an elliptic curve has 2-torsion points. If the elliptic
curve has equation $y^2=x^3+Ax+B$, the abscissas of the points of
order 2, are given by the solutions of te cubic $x^3+Ax+B$ over the
finite field. A study of elliptic curves with 2-torsion points can
be found in
\cite{IMJ} using elliptic curves in Legendre form.\\

For an integer $n \neq 2$, the n-torsion points of an elliptic curve
have as abscissas the roots of the classical division polynomials.
These can be constructed recurrently:
$$
\begin{array}{ll}
P_{-1}=-1, \quad P_{0}=0, \quad  P_{1}=1,\quad P_{2}(x,y)=2y\\
P_{3}=3x^4+6ax^2+12bx-a^2\\

%\item $P_{4}=4y(x^6+5ax^4+20bx^3-5a^2x^2-4abx-8b^2-a^3)$

P_{2k}=\frac{P_{k}}{2y}(P_{k+2}P_{k-1}^2-P_{k-2}P_{k+1}^2),
P_{2k+1}=P_{k+2}P_{k}^3-P_{k+1}^3P_{k-1}
\end{array}
$$

But, with these polynomials it is not possible to characterize all
elliptic curves with a point of order $n$. In \cite{mestre}, Mestre
gives equations for $\_0(N)$, $N=2,3,5,7, 13$ which parametrizes the
pairs $(E,C)$ of an elliptic curve $E$ and a subgroup $C$ of order
$N$. He gives equations of the modular curve via the j- invariant of
the curve. Our purpose is to present a family of the elliptic curves
with a rational subgroup of order 3, but in their Weiertrass normal
form. The interest of these curves comes from the study of volcanoes
of 3- isogenies (see \cite{volcanes}).

The structure of the paper is the following. In section
\ref{3racional} we present the elliptic curves with, at least, a
point of order 3. As it is known (\cite{Hus}), an elliptic curve
with a 3-torsion point over a finite field, admits a model of the
form $ y^2+a_1x+a_3y=x^3$. Nevertheless, we can consider an elliptic
curve isomorphic to the former one that depends only on one
parameter (after fixing a non cube). We distinguish between the
cases when the 3-torsion subgroup is cyclic or not. After this
classification, we determine the number of isomorphism classes of
elliptic curves with this property.

Section  \ref{nine} is dedicated to study when a rational point $P$
of order 3 on an elliptic curve (among those mentioned in the
previous section) has a point $Q$ which \emph{trisects} it (i.e.
$3Q=P$). This follows easily from \cite{MMRV} and \cite{Jac}.

Finally, in section \ref{nocyclic} we determine a family of
representatives of all isomorphic classes of elliptic curves with a
rational subgroup of order 3 whose points are not rational. For this
purpose, we compare the 3 division polinomial of a curve with that
of its twisted curve, to conclude that both factorize in the same
way. Thus, the curves in section \ref{3racional} are in
correspondence with the desired curves in this section.

\section{Elliptic curves with cardinal a multiple of 3}\label{3racional}

 Let $E$ be an elliptic curve over a finite field
${\mathbb{F}_q} $. It's  well known that $E$ admits a model in
Weierstrass normal form (\cite{Sil}):

$$
E \, :\, y^2+a_1xy+a_3y=x^3+a_2x^2+a_4x+a_6
$$

The rational points of order 3 can be obtained from the 3-division
polynomial. The roots of this polynomial correspond to the abscissas
of the 3-order points \footnote{Although these abscissas are defined
over  ${\mathbb{F}_q} $, it is possible that the point is not. So,
they are points over ${\mathbb{F}}_{q^2} $}. Let $E$ be an elliptic
curve  with, at least, a point  of order 3, we can suppose, without
loss of generality, that this point is $(0,0) $. In this situation,
$E $ admits an equation (translating the 3-order point to the
origin) of the type
\begin{equation}\label{punt3}
  y^2+a_1x+a_3y=x^3
\end{equation}

%and the tangent to the curve at $(0,0) $ is $y=0$. \\

The discriminant of an elliptic curve of equation (\ref{punt3}) is
$\triangle=(a_1^3-27a_3)a_3^3 $ and the points $(0,0) $ and
$(0,-a_3) $ are of order $3 $.\\

Equation (\ref{punt3}) can be simplified  as follow:

\begin{lema}\label{curvasisomorfascubos}

Let $E $ be an elliptic curve  over ${\mathbb{F}_q} $ with
equation $y^2+a_1x+a_3y=x^3$ and let $b_0 \in {\mathbb{F}_q} $ be
a fixed non cube. If $b $ is a cube, then $E $ is isomorphic to a
curve of type $y^2+axy+y=x^3 $. Otherwise,  then $E $ is
isomorphic to $y^2+axy+b_0y=x^3 $ or $y^2+axy+b_0^2y=x^3 $.

Moreover, if $q \equiv -1 \bmod{3} $, all elliptic curves with a
point of order 3 admit an equation of the form $y^2+axy+y=x^3 $.
\end{lema}

\demo \\

Let us suppose that $b $ is a cube. Let $E$ be the elliptic curve
with equation $y^2+axy+\alpha^3 y=x^3 $. Replacing $(x,y)$ by
$(\alpha^2x',\alpha^3 y') $, it gives an equation
of the form $y^2+(a/\alpha) xy +y=x^3 $. \\

If $b $ is not a cube, let $\rho $ be a 3-order primitive root of
the unity and let $\chi $ be a cubic character over ${\mathbb{F}_q}
$. It is clear that $\chi(b) = \rho $ or $\rho^2 $. We can suppose
that $\chi(b_0)=\rho $ ( otherwise, we take its square). If
$\chi(b)=\rho $, replacing $(x,y)$ by $(u^2x',u^3 y') $,
$u=\sqrt[3]{b/b_0} $; the elliptic curve has equation
$y^2+axy+b_0y=x^3 $. If $\chi(b)=\rho^2 $, replace $(x,y)$ by
$(u^2x',u^3 y')$, $u=\sqrt[3]{b/b_0^2} $ to get the
desired equation. \\

Finally, over ${\mathbb{F}_q} $, $q \equiv -1 \bmod{3} $, all
elements are cubes, therefore, all elliptic curves with a 3-order
point,  admit an equation of the type $y^2+axy+y=x^3 $. \finde

\subsection{Elliptic curves with cyclic 3-torsion subgroup}

The 3-division polynomial of  an elliptic curve  with equation
$y^2+axy+y=x^3 $ is $P_3(x)=x(3x^3+a^2x^2+3ax+3) $. If $q \equiv -1
\bmod{3} $, this polynomial has two roots over ${\mathbb{F}_q} $.
Nevertheless, only the root 0 corresponds to an abscissa of a
3-order point. The other root corresponds to a 3-order point defined
over a quadratic extension of $\mathbb{F}_q$. This means that there
exists a 3-order rational subgroup with no rational points. So, we
have the following result:

\begin{Th}\label{card3q2}
Let ${\mathbb{F}_q} $, $q \equiv -1 \bmod{3} $ a finite field, then
the only elliptic curves with 3-order rational points are isomorphic
to a curve with equation $y^2+axy+y=x^3 $, $a\neq 3 $. Moreover, the
3-torsion group is, in this case, cyclic \footnote{For all these
elliptic curves, the 3-order points are $(0,0) $ and $- (0,0)=(0,-1)
$}.
\end{Th}

Taking into account these equations for the elliptic curves, we can
count the number of isomorphism classes of elliptic curves over
${\mathbb{F}_q}$, $q \equiv -1 \bmod{3}$, with a point of order
three.

\begin{prop}
Let ${\mathbb{F}_q} $, with $q \equiv -1 \bmod{3} $ a finite field,,
there exist $q-1$ isomorphic classes of elliptic curves with
cardinal a multiple of 3. A parametrization of these classes is
given by
\begin{equation}
  \{y^2+axy+y=x^3 \, | \, a \neq 3 \}
\end{equation}
\end{prop}

\demo \\

By theorem \ref{card3q2}, let $E\,:\, y^2+axy+y=x^3 $ be an elliptic
curve over ${\mathbb{F}_q} $, $q \equiv -1 \bmod{3}$, with a
3-torsion point. The Weierstrass normal form of this curve is
$y^2=x^3+Ax+B $ where

\begin{equation} A_a=-\frac{a(-24+a^3)} {48} \qquad
B_a=\frac{1}{4}+\frac{a^6}{864}-\frac{a^3}{24}.
\end{equation}

Over  $q \equiv -1 \bmod{3} $, every element is a  cube, and for
each one, one only cubic root exists. Therefore, we can only
consider the curves $y^2+\sqrt[3]{a}xy+y=x^3 $ (this equation makes
the calculations easier). The j-invariant of this curve is
$-16 \frac{a(-24+a)^3}{a-27} $.\\

Let $y^2+\sqrt[3]{b}xy+y=x^3 $ be another elliptic curve of this
type. A necessary condition for the elliptic curve to be isomorphic
to the previous one  is that both possess the same j-invariant. So
$b $ will be a root of
$$
X^3+(a-72)X^2+(a^2-72a+1728)X-27\frac{(-24+a^3)} {to-27}.
$$

The only root of this polynomial over the finite field considered
is
\begin{equation}\label{biso}
b=\frac{72-a}{3}-\frac{(-36+a)\sqrt[3]{a}} {3\sqrt[3]{-27+a}}
+\frac{2}{3}\sqrt[3]{-27+a}\sqrt[3]{a^2}. \end{equation}

Also, if they are isomorphic, there exists $u \in {\mathbb{F}_q} ^*$
such that $A_a=u^4A_b $, $B_a=u^6B_b $. Replacing (\ref{biso}) in
this second equation, one has that
$$
u^6=-3\frac{(\sqrt[3]{-27+a}) ^2}{(\sqrt[3]{-27+a}-\sqrt[3]{a}) ^2}.
$$

But, in this situation, $-3 $  is not a square and thus, isomorphic
elliptic curves to $y^2+\sqrt[3]{a}xy+y=x^3 $ do not exist.
 \finde

For finite fields ${\mathbb{F}_q}$ with $q \equiv 1 \bmod{3}$, we
have to consider all the different types of elliptic curves shown in
lemma \ref{curvasisomorfascubos}:

\begin{equation}\label{curva1}
\{E^i_a\,:\,y^2+axy+b_0^iy=x^3 \,; \, i=0,1,2 \}. \end{equation}

where $\chi(b_0)=\rho $ ($\chi $ is a cubic character
over ${\mathbb{F}_q} $ and $\rho $ a cubic root of the unity).\\

From the 3-division polynomial, it is easy to see that these curves
will have cyclic 3-torsion subgroup if
and only if $\chi(27b_0^i-a^3)\neq 1 $. \\

By Lemma \ref{curvasisomorfascubos}, if $E^i_{a_1} \cong E^j_{a_2} $
then  $i=j $. We also have:

\begin{lema}\label{coleccisomorfas}
Let $E^i_{a_1}, E^i_{a_2} $ be two elliptic curves over
${\mathbb{F}_q} $, $q \equiv 1 \bmod{3} $ such that
$\chi(27b_0^i-a_1^3), \chi(27b_0^i-a_2^3) \neq 1 $ (both curves
possess cyclic 3-torsion subgroup). $E^i_{a_1} \cong E^i_{a_2} $
if and only if:
$$
a_2 \in \{a_1,\rho a_1, (\rho+1)a_1 \}
$$
\end{lema}

\demo \\

It follows from the j-invariants of both curves. These are the only
possibilities over ${\mathbb{F}_q} $, $q \equiv 1 \bmod{3} $ and it
is easy to prove that these elliptic curves are isomorphic. \finde

As $\chi(27b_0^i-a^3) \neq 1 $,  determining the number of possible
values for $a $,  is equivalent to determine how many values of $a$
exist with $a^3-b_0^i $ been a cube. The following result will be
essential.

\begin{lema}
Let ${\mathbb{F}_q} $ be a finite field, with $q \equiv 1 \bmod{3} $
and let $A \equiv 1 \bmod{3} $ be such that $4q=A^2+27B^2 $ for some
$B $ (with these conditions $A $  is unique).

\begin{enumerate}
\item The number of solutions of $x^3+y^3=1 $  is exactly $q-2+A$.

\item The number of pairs $(x,y) $ such that $x^3+y^3 $  is not a
cube t is $\frac{(q-1)(2q-4-A)}{3} $.

\end{enumerate}
\end{lema}
\demo \\

A proof of part i) can be found in \cite{IrRo}.\\

 For each $a$, the
number of solutions of $x^3+y^3=a^3 $ is also $q-2+A $, except for
$a=0 $, in that case there are $3q-2$.  Therefore, since there exist
$(p-1)/3$ cubes different from 0, there will be $q^2-(\frac{(q-1)}
{3}(q-2+A)+3q-2)=\frac{(q-1)(2q-4-A)} {3} $ pairs $(x,y) $ such that
$x^3+y^3 $ is not a cube. \finde

Using this result we can compute the number of isomorphism classes
of elliptic curves over ${\mathbb{F}_q} $, $q \equiv 1 \bmod{3} $
with cyclic 3-torsion.

\begin{prop}\label{nociclic}
Let ${\mathbb{F}_q} $, $q \equiv 1 \bmod{3} $ a finite field. There
are $(2q+4)/3$ isomorphic classes of elliptic curves with cyclic
3-torsion subgroup. Let $b_0 \in {\mathbb{F}_q}$ be a non cube. A
family of representatives  is given by:
\begin{equation}\label{toda3torsionrac} \{y^2+m_axy+b_0^iy=x^3
\,; \, i=0,1,2, \chi(a^3-27b_0^i)\neq 1 \} \end{equation} with
$m_a=min\{a,\rho a, (\rho+1)a\}$.
\end{prop}

\demo \\

First, we compute the number of isomorphism classes of elliptic
curves of the form $y^2+axy+y=x^3 $. This curve has a cyclic
3-torsion if and only if $a^3-27 $ is not a cube. By lemma
\ref{coleccisomorfas}, the curves $y^2+\rho a xy+y=x^3 $ and
$y^2+\rho^2 axy+y=x^3 $ are isomorphic to the previous one.
Therefore, the number of isomorphism classes of these curves is
given by the number of cubes $a^3 $, such that
$\chi(a^3-27)=\chi(a^3-1)\neq 1 $.\\

By the previous lemma, there exist $q-2+A $ pairs $(x,y) $ such that
$x^3+y^3=1$. To study the number of elements $x$ such that both $x $
and $x-1 $ are cubes, it is necessary to keep in mind the following
facts. The pairs $(0,1),(0, \rho) $ and $(0,\rho^2) $ produce the
same cube $0 $. The same fact is true  for $(1,0),(\rho,0) $ and
$(\rho^2,1) $ and finally for $(\rho^i x,\rho^j y) $, $0\leq
i,j \leq 2 $.\\

Therefore, there exist $\frac{q+10+A}{9} $ cubes $x^3 $, such that
$x^3-1 $ is a cube. In consequence, there are  $(2q-4+A)/9 $
cubes $x^3 $ such that $x^3-1$ is not.\\

In a similar way, it can be shown that there exist $(4q+16+A)/9 $
isomorphism classes of elliptic curves of the type
$y^2+a^3xy+b_0^iy=x^3 \,; \,$ with $ i=1,2,$ and $
\chi(a^3-b_0^i)\neq 1 $ and the result follows after adding this
number to the previous one. \finde

\subsection{Elliptic curves with non cyclic 3-torsion subgroup}

As shown in the previous section, a necessary condition for an
elliptic curve to have non cyclic 3-torsion subgroup is that it be
defined over ${\mathbb{F}_q} $ with $q \equiv 1 \bmod{3}$. If $E$ is
defined by the equation $ y^2+axy+by=x^3$, it is also necessary that
$\chi(a^3-27b) \neq 1$. We can change the equation to distinguish
between those with non cyclic subgroup as follows:

\begin{prop}
The family of elliptic curves
$$
E_a\,:\,y^2+(3a-1)xy+a(\rho-1)(a-\frac{\rho+1}{3}) y=x^3
$$
such that $a(a-(\rho+1)/3)(a-\rho/3)\neq 0 $ corresponds to all
the  elliptic curves defined over ${\mathbb{F}_q} $ with non
cyclic 3-torsion subgroup. \end{prop}

\demo \\

Let $y^2+axy+by=x^3 $ be an elliptic curve over ${\mathbb{F}_q} $
with non cyclic 3-torsion subgroup. It is clear that (by
construction) one of the points of order 3 is the origin. Let
$(x_0,y_0) $ be another point of order 3 which does not belong to
$<(0,0)> $ (that is $x_0\neq 0 $). Let $y=\lambda x +\mu $ be the
tangent line to the curve in $(x_0,y_0) $ ($\lambda \neq 0 $ because
this point
is not the origin).\\

The change of variables $(x,y) \rightarrow (\lambda^2 x, \lambda^3
y) $  transforms the equation of the elliptic curve  in to another
of the same type.  The  point of order 3 is now $(x_0/\lambda^2,
y_0/\lambda^3) $ and the tangent line at this point has slope
$1$.\\

The result follows easily by considering the conditions in the cubic
such that both $y=0$ and $y=x+u$ have a triple intersection point
with the cubic (this can be found in \cite{Hus}). \finde

Let us see now when two elliptic curves $E_a $ and $E_b $ are
isomorphic. It is clear that $a \not \in \{0, \rho/3, (\rho+1)/3 \}
$, otherwise it would be a singular cubic. Using the expressions for
the j-invariants and the equations in  Weierstrass normal form of
these curves we can prove that:

\begin{lema}
Let ${\mathbb{F}_q} $ be a finite field with $q \equiv 1 \bmod{3}
$ and let $E_a $ and $E_b $ be two elliptic curves over
${\mathbb{F}_q} $. The equations for $E_a $ and $E_b $  are
\begin{eqnarray} y^2+(3a-1)xy+a(\rho-1)(a-\frac{\rho+1}{3}) y=x^3
\nonumber \\y^2+(3b-1)xy+b(\rho-1)(b-\frac{\rho+1}{3})
y=x^3\nonumber \end{eqnarray}
 respectively
(such that $a,b \not \in \{0, \rho/3, (\rho+1)/3 \} $). $E_a \cong
E_b $ if and only if
$$
\begin{array}{lr}
b \in G_a=\{a, \frac{a(1+\rho)} {(3a-\rho)},
\frac{a\rho}{(3a-1-\rho)},
\frac{-1}{9a},\frac{-(1+\rho)(3a-1-\rho)} {3},
\frac{\rho(3a-\rho)} {9a}, \frac{\rho(3a-1-\rho)} {3(3a-\rho)},
\\\hspace * {.4cm} & \\\hspace*{2cm}\frac{-\rho}{3(3a-1-\rho)}, \frac{(1+\rho)(3a-\rho)}
{3(3a-1-\rho)}, \frac{(1+\rho)(3a-1-\rho)} {9a},\frac{(1+\rho)}
{3(3a-\rho)}, \frac{\rho*(3a-\rho)} {3} \}. \end{array}
$$
\end{lema}

$G_a $ has a group structure and it is isomorphic to
${\mathbb{Z}}/2{\mathbb{Z}} \times {\mathbb{Z}}/6{\mathbb{Z}} $.
$G_a $ acts over ${\mathbb{F}_q} - \{0, \rho/3, (\rho+1)/3 \} $.
Then, the number of isomorphism classes of elliptic curves with non
cyclic  coincides whit the number of orbits under this action. To
compute this number we  use  the well-known Burnside formula

$$
\sharp{\mbox { of orbits}} = \frac{1}{|G_a |} \sum_{g \in G_a}|N^g |
$$
where $N^g $ is the set of elements fixed by $g $. We have:
\begin{center}
\begin{tabular}{||l|c||l|c ||}
  \hline
  \hline
   $g $& $N^g $& $g$&$N^g $\\  \hline
$g$&${\mathbb{F}_q} - \{0, \rho/3, (\rho+1)/3 \} $&
$\frac{\rho(3a-1-\rho)} {3(3a-\rho)} $& $\frac{\rho+\sqrt{-\rho}}
{3}, \frac{\rho-\sqrt{-\rho}} {3} $\\\hline $\frac{a(1+\rho)}
{(3a-\rho)} $& $\frac{1+2\rho}{3} $& $\frac{-\rho}{3(3a-1-\rho)}
$& $\frac{1}{3} $\\\hline $\frac{a\rho}{(3a-1-\rho)} $&
$\frac{1+2\rho}{3} $& $\frac{(1+\rho)(3a-\rho)} {3(3a-1-\rho)} $&
$\frac{\rho+1+\rho\sqrt{-1}} {3}, \frac{\rho+1-\rho\sqrt{-1}} {3}
$\\\hline $\frac{-1}{9a}$&$\frac{\sqrt{-1}} {3}, \frac{-\sqrt{-1}}
{3} $& $\frac{(1+\rho)(3a-1-\rho)} {9a}$&$\frac{1}{3} $\\\hline
$\frac{-(1+\rho)(3a-1-\rho)} {3}$&$\frac{1}{3} $& $\frac{(1+\rho)}
{3(3a-\rho)} $& $- \frac{1}{3} $\\\hline $\frac{\rho(3a-\rho)}
{9a}$&$-\frac{1}{3} $& $\frac{\rho*(3a-\rho)}
{3}$&$\frac{1+2\rho}{9} $\\\hline \hline \end{tabular}
\end{center}

Considering whether these elements belong to ${\mathbb{F}_q} - \{0,
\frac{\rho}{3},\frac{(\rho+1)} {3} \} $ and the Burnside formula, we
have the following result.

\begin{prop}

Let ${\mathbb{F}_q} $ be a finite field with $q \equiv 1 \bmod{3}
$. Over ${\mathbb{F}_q} $ there exist $\frac{q+12-(q \bmod{12})}
{12} $ isomorphic classes of elliptic curves with non cyclic
3-torsion subgroup.

\end{prop}

\section{Elliptic curves with rational subgroups of order 3
containing no rational points}\label{nocyclic}

In this section we are interested in elliptic curves that have no
rational points of order three, but they have a rational subgroup of
order 3. These curves appear when we consider rational isogenies of
degree three. This is because the points of this rational subgroup
are not invariant under the action of the Frobenius
endomorphism, whereas the subgroup  is.\\

Let $E \, :\, y^2=x^3+Ax+B $ be an elliptic curve defined over a
finite field ${\mathbb{F}_q} $ and let $G \subset E({\mathbb{F}_q})
$ be a subgroup of order $3$ whose points are not rational (except
the point of the infinity). Let $\sigma: {\mathbb{F}_q} \rightarrow
{\mathbb{F}_q} $, be the Frobenius endomorphism: $\sigma((x,y)) =
(x^q,y^q) $. If $P \in G $, then $\sigma(P)=P $or $- P $. If it
coincides with $P $, it is clear that the point $P $ is rational,
which is absurd. Therefore, $\sigma(P)=-P $ and hence
$$
(x^q,y^q)=(x,-y).
$$

From this, we deduce that the abscissa of  $P $  should be
rational and in consequence the ordinate should be defined over
$\mathbb{F}_{q^2} $. Then $\sharp(E(\mathbb{F}_{q^2}))
$ is a multiple of 3. \\

By Weil theorem, if $\sharp(E(\mathbb{F}_{q})) =
m({\mathbb{F}_q})=q+1-t $, then $\sharp(E(\mathbb{F}_{q^2})) =
m(\mathbb{F}_{q^2}) = q^2+1-t^2+2q $. If $q \equiv -1 \bmod{3} $,
$m(\mathbb{F}_{q^2}) $ is a multiple of 3 if and only if
$m({\mathbb{F}_q}) $  is . Therefore, all the elliptic curves with
rational subgroups of order 3 over ${\mathbb{F}_q} $, $q \equiv -1
\bmod{3} $, have been already studied in section
\ref{3racional}.\\

Otherwise, if $q \equiv 1 \bmod{3} $, $m(\mathbb{F}_{q^2}) $ is a
multiple of 3 if and only if $t \equiv 1, 2 \bmod{3}$. In the first
case, the  curve defined over ${\mathbb{F}_q} $  has cardinal
multiple of 3 (already studied in section \ref{3racional}). In the
second case, $m({\mathbb{F}_q}) $ is not a multiple of 3. In
consequence, we must only study elliptic curves with points of order
3 defined over $\mathbb{F}_{q^2} $ but not over ${\mathbb{F}_q} $.
Moreover, this implies that there exist a one-to-one correspondence
between curves with rational points of order 3 and those with
3-order subgroups with no rational points (ones are
the twisted of the others). \\

The 3-division polynomial (the one whose roots are the abscissas of
the  points of order three) of $E $  is

$$
P_{A,B}=x^4+2Ax^2+4Bx-\frac{A^2}{3}.
$$

The following theorem gives us some information about  the
factorization of $P_{A,B}$ over ${\mathbb{F}_q} $, $q \equiv 1
\bmod{3} $.

\begin{Th}[\cite{Leo}]
Let $p$ be an odd prime and let $f(x) $ be a monic polynomial of
degree $n $ and  discriminant $D $ without multiple roots. Let
$f(x)=f_1(x)f_2(x)\cdots f_r(x) \bmod{p} $ be the factorization of
$f$ over ${\mathbb{F}}_p$. Then $n \equiv r \bmod{2} $ if and only
if $(D/p)=1 $.
\end{Th}

The discriminant of $P_{A,B} $  is $-2^8(4A^3+27B^2)^2/27 $. Over
${\mathbb{F}_q} $, $q \equiv 1 \bmod{3} $, $-1/3 $ is a square. So,
one has that $P_{A,B} $ factors  into two polynomials or it has all
its roots in ${\mathbb{F}_q} $. If it factors into two polynomials,
the possibilities for its degree are $(2,2)$ or $(3,1)$. In the
latter case, $P_{A,B} $ has a root in ${\mathbb{F}_q} $. \\

Using the results of \cite{Sko}, these three types  can be
distinguished in the following way. If $\chi(2(27B^2+4A^3)) \neq 1
$, then $P_{A,B}=(x-x_0)f_2(x) $ with $f_2(x) $ an irreducible
polynomial of degree 3. Otherwise, we define
$y_0=2\sqrt[3]{2(27B^2+4A^3)} /3 $. If $y_0+16A/3 $ and $\rho
y_0+16A/3 $ are both squares, then $P_{A,B} $ factors completely
over ${\mathbb{F}_q} $. Otherwise, $P_{A,B} $ factors in two
irreducible
quadratic polynomials.\\

With these results, we have that the type ot factorization of
$P_{A,B} $ coincides with the that of $P_{t^2A,t^3B} $, where $t $
is a non quadratic residue.  That corresponds to the 3-division
polynomial of the twisted curve. \\

As we have shown previously,  elliptic curves with rational points
of order three are the twisted of those with subgroups of order 3
with no rational points. Thus, if $y^2=x^3+Ax+B $ is an equation for
an elliptic curve with 3-torsion points, then $y^2=x^3+t^2Ax+t^3B $
($t $ a non quadratic residue), is a curve with a rational subgroup
of order 3 with no rational points. As we have studied in section
\ref{3racional} the curves of the first type,  determine those of
the second type from their Weierstrass form.

\begin{prop}
Let ${\mathbb{F}_q} $, $q \equiv 1 \bmod{3} $ be a finite field and
$b_0 \in {\mathbb{F}_q} $, such that $\chi(b_0)=\rho $. Let $t $ be
a non quadratic residue. There exist $(2q+4)/3 $ isomorphism classes
of elliptic curves over ${\mathbb{F}_q} $ with only one rational
subgroup of order three, whose points are defined over
$\mathbb{F}_{q^2}\setminus{\mathbb{F}_q} $. A family of
representatives is given for $y^2=x^3+t^2Ax+t^3B $ with

$$
A_a=-\frac{a^3(a^9-24b_0^i)} {48} \qquad B_a=\frac{b_0^{2i}}
{4}+\frac{a^{18}} {864}-\frac{a^9b_0^i}{24}.
$$
$$
\quad i=0,1,2 \quad \chi(a^3-27b_0^i)\neq 1
$$
\end{prop}

\demo \\

Elliptic curves of equation (\ref{toda3torsionrac})
$\{y^2+a^3xy+b_0^iy=x^3 \,; \, i=0,1,2, \chi(a^3-27b_0^i)\neq 1 \}
$ correspond to those curves with a non cyclic 3-torsion subgroup
whose points are defined over ${\mathbb{F}_q} $. Therefore, their
twisted curves will have a rational 3-torsion subgroup of points
defined in a quadratic extension. \\

The Weierstrass normal form of elliptic curves with equation
(\ref{toda3torsionrac})  is
$$
y^2=x^3+ \biggl(- \frac{a^3(a^9-24b_0^i)} {48}\biggr) x+
\frac{b_0^{2i}} {4}+ \frac{a^{18}} {864}-\frac{a^9b_0^i}{24}.
$$

In consequence, the equation of the twisted curve coincides with the
desired one.\\

Let  $y^2=x^3+t^2Ax+t^3B $ and $y^2=x^3+t^2A'x+t^3B' $ be two of
these curves. They will be isomorphic over ${\mathbb{F}_q} $, if and
only if there exists $u \in {\mathbb{F}_q}^* $ such that $u^4=A/A',
\, u^6=B/B' $. Or equivalently, if and only if the curves
$y^2+a^3xy+b_0^iy=x^3, \, y^2+a'^3xy+b_0^iy=x^3 $ are isomorphic.
The number of isomorphism classes follows from proposition
\ref{nociclic}.\finde

Two curves of this type will be isomorphic if and only if
$$
a_2 \in\{a_1,\rho a_1, (\rho+1)a_1 \}.
$$

We can prove, in the same way, the following result:

\begin{prop}
Let ${\mathbb{F}_q} $ be a finite field with $q \equiv 1 \bmod{3} $.
Over ${\mathbb{F}_q} $ there exist $\frac{q+12-(q \bmod{12})} {12} $
isomorphism classes of elliptic curves  with 4 rational subgroups of
order 3 composed by points defined over $\mathbb{F}_{q^2} $. A
collection of representatives is given by $y^2=x^3+t^2Ax+t^3B $ with
$t $ a non quadratic residue and
$$
\begin{array}{l}
A=- \frac{(9a-1-2\rho)(3a-1-2\rho)(3a-1)(3a+1)} {144}
\\B=\frac{(1+9a^2)(9a^2-6a-6\rho to-1)(9a^2-6\rho to-1)} {864} \end{array}
$$

\end{prop}

Two curves of this type are isomorphic (the same as it happened
with those with 8 rational points of order 3) if and only if
$$
\begin{array}{lr}
b \in G_a=\{a, \frac{a(1+\rho)} {(3a-\rho)},
\frac{a\rho}{(3a-1-\rho)},
\frac{-1}{9a},\frac{-(1+\rho)(3a-1-\rho)} {3},
\frac{\rho(3a-\rho)} {9a}, \frac{\rho(3a-1-\rho)} {3(3a-\rho)},
\\\hspace * {.4cm} & \\\hspace*{2cm}\frac{-\rho}{3(3a-1-\rho)}, \frac{(1+\rho)(3a-\rho)}
{3(3a-1-\rho)}, \frac{(1+\rho)(3a-1-\rho)} {9a},\frac{(1+\rho)}
{3(3a-\rho)}, \frac{\rho*(3a-\rho)} {3} \}. \end{array}
$$

\end{document}